\documentstyle[12pt]{article}
\newtheorem{th}{\indent Theorem}
\newcommand{\R}{{\bf R}}
\newcommand{\M}{{\cal M}}
\begin{document}
\begin{center}
{\large\bf
A new combinatorial characterization of the minimal cardinality
of a subset of $\R$ which is not of first category}
\par\vspace{.3cm}\par
{\bf Apoloniusz Tyszka}
\end{center}
\begin{abstract}
\def\thefootnote{}
Let $\M$ denote the ideal of first category subsets of $\R$.
We prove that $\min \{\rm card \ X: X\subseteq\R, X\not\in \M\}$
is the smallest cardinality
of a family $S\subseteq\{0,1\}^\omega$ with the property that
for each
$f:\omega\longrightarrow\bigcup_{n\in\omega}\{0,1\}^n$ there exists a
sequence $\{a_n\}_{n\in\omega}$ belonging to $S$ such
that for infinitely many $i\in\omega$ the infinite sequence
$\{a_{i+n}\}_{n\in\omega}$ extends the finite sequence $f(i)$.
\\
\\
We inform that $S \subseteq \{0,1\}^\omega$ is not of first category
if and only if for each
$f:\omega\longrightarrow\bigcup_{n\in\omega}\{0,1\}^n$ there exists a
sequence $\{a_n\}_{n\in\omega}$ belonging to $S$ such
that for infinitely many $i\in\omega$ the infinite sequence
$\{a_{i+n}\}_{n\in\omega}$ extends the finite sequence $f(i)$.
\end{abstract}

Let $\M$ denote the ideal of first category subsets of $\R$. Let
$\M(\{0,1\}^\omega)$ denote
the ideal of first category subsets of the Cantor space
$\{0,1\}^\omega$. Obviously:
$$
\arraycolsep 2pt
\begin{array}{lcl}
\mbox{non}(\M):&=&\min\{\mbox{card } X: X\subseteq\R, X\not\in \M\}\\
&=& \min\{\mbox{card } X: X\subseteq \{0,1\}^\omega,
X\not\in\M(\{0,1\}^\omega)\}\\
\end{array}
\leqno(\ast)
$$
Let $\forall^\infty$ abbreviate "for all except finitely many". It is
known (see [1], [2] and also [3]) that:
$$
\begin{array}{c}
\mbox{non}(\M) = \\ [.2cm]
\min\{\mbox{card } F: F\subseteq\omega^\omega\mbox{ and } \neg\ \exists
g\in\omega^\omega \ \forall f\in F \ \forall^\infty k \ g(k)\neq f(k)\}
\end{array}
$$
For $S\subseteq\{0,1\}^\omega$ we define the following property
$(\ast\ast)$:
\begin{center}
\begin{tabular}{ll}
&for each $f:\omega\rightarrow\bigcup_{n\in\omega}\{0,1\}^n$ there exists
a sequence $\{a_n\}_{n\in\omega}$ \\
$(\ast\ast)$&belonging to $S$ such that for infinitely
many $i\in\omega$ the \\
&infinite sequence $\{a_{i+n}\}_{n\in\omega}$ extends the
finite sequence $f(i)$.
\end{tabular}
\end{center}
\begin{th}
If $S\subseteq\{0,1\}^\omega$ is not of first category then $S$ has the
property $(\ast\ast)$.
\end{th}
{\it Proof.} Let us fix
$f:\omega\longrightarrow\bigcup_{n\in\omega}\{0,1\}^n$. Let $S_k (f)$
$(k\in\omega)$ denote the set of all sequences $\{a_n\}_{n\in\omega}$
belonging to $\{0,1\}^\omega$ with the property that there exists
$i\geq k$ such that the infinite sequence $\{a_{i+n}\}_{n\in\omega}$
extends the finite sequence $f(i)$.
\par
Sets $S_k (f)$ $(k\in\omega)$ are open and dense. In virtue of the Baire
category theorem $\bigcap_{k\in\omega} S_k (f)\cap S$ is non-empty i.e.
there exists a sequence $\{a_n\}_{n\in\omega}$ belonging to $S$
such that for infinitely many $i\in\omega$ the infinite sequence
$\{a_{i+n}\}_{n\in\omega}$ extends the finite sequence $f(i)$. This
completes the proof.
\\
\\
{\bf Note.} The author recently proved that if
$S \subseteq \{0,1\}^\omega$ has the property ($\ast\ast$) then $S$
is not of first category; the proof will appear in a separate
preprint. From this result and Theorem 1 we obtain the
following characterization:
\\
$S \subseteq \{0,1\}^\omega$ is not of first category if and only if
$S$ has the property ($\ast\ast$).
\\
Let us note that from the above characterization we may deduce all next
results; therefore all next proofs are unnecessary.
\\
\\
\begin{th}
If $S\subseteq\{0,1\}^\omega$ has the property $(\ast\ast)$ then
{\rm card} $S\geq {\rm non}(\M)$.
\end{th}
{\it Proof.} For a sequence $\{a_n\}_{n\in\omega}$ belonging
to $S$ we define:
$$
\tilde{a}_n := \sum_{i=n}^\infty\frac{a_i}{2^{i-n+1}}\in [0,1].
$$
The following Observation is easy.
\par\vspace{.3cm}\par
{\bf Observation.} Assume that $S\subseteq\{0,1\}^\omega$ has the
property $(\ast\ast)$. We claim that for each sequence
$\{U_k\}_{k\in\omega}$ of
non-empty open sets satisfying $U_k \subseteq (0,1)$ $(k\in\omega)$ there
exists a sequence $\{a_k\}_{k\in\omega}$ belonging to $S$ such that for
infinitely many $k\in\omega$ $\tilde{a}_k \in U_k$.
\vskip 0.5truecm
\par
There exists a sequence $\{(c_i,d_i)\}_{i\in\omega}$ of non-empty
pairwise disjoint intervals satisfying
$$
\bigcup_{i\in\omega} (c_i,d_i)\subseteq (0,1).
$$
We assign to each $\{a_n\}\in S$ the function
$s_{\{a_n\}}:\omega\longrightarrow\omega$ according to the following rules
(cf.[4]):
\begin{itemize}
\item[1)] if $\tilde{a}_k \not\in\bigcup_{i\in\omega} (c_i,d_i)$ then
$s_{\{a_n\}}(k)=0$,
\item[2)] if $\tilde{a}_k \in\bigcup_{i\in\omega} (c_i,d_i)$ then
$s_{\{a_n\}}(k)$ is the unique $i\in\omega$ such that $\tilde{a}_k \in
(c_i,d_i)$.
\end{itemize}
Suppose, contrary to our claim, that card $S < \mbox{ non}(\M)$. It
implies that the cardinality of the family $\{s_{\{a_n\}}:\{a_n\}\in
S\}\subseteq\omega^\omega$ is also less than
$$
\begin{array}{c}
\mbox{non}(M) = \\ [.2cm]
\min \{\mbox{card } F: F\subseteq\omega^\omega \mbox{ and } \neg \ \exists
\ g\in\omega^\omega \ \forall f\in F \ \forall^\infty k \ g(k)\neq f(k)\}
\end{array}
$$
Therefore, there exists a function $g:\omega\longrightarrow\omega$ such
that for each sequence $\{a_n\}\in S \ \forall^\infty k \ g(k)\neq
s_{\{a_n\}}(k)$. We define $U_k := (c_{g(k)},d_{g(k)})\subseteq (0,1) \
 (k\in\omega)$. If $\{a_n\}\in S$ then the set
$A_{\{a_n\}}:=\{k\in\omega: g(k)=s_{\{a_n\}}(k)\}$ is finite and for each
$k\in\omega\setminus A_{\{a_n\}} \ \tilde{a}_k \not\in U_k$. It
contradicts the thesis of the Observation which ensures that there
exists a sequence $\{a_k\}_{k\in\omega}$ belonging to $S$ such that for
infinitely many $k\in\omega \ \tilde{a}_k \in U_k$.
This completes the proof of Theorem 2.
\par\vspace{.3cm}\par
{\bf Corollary.} From $(\ast)$, Theorem 1 and Theorem 2 follows that
non$(\M)$ is the smallest cardinality of a family
$S\subseteq\{0,1\}^\omega$ with the property that for each
$f:\omega\longrightarrow\bigcup_{n\in\omega}\{0,1\}^n$ there exists a
sequence
$\{a_n\}_{n\in\omega}$ belonging to $S$ such that for infinitely many
$i\in\omega$ the infinite sequence $\{a_{i+n}\}_{n\in\omega}$ extends
the finite sequence $f(i)$.
\par
{\bf Remark.} Another (not purely combinatorial) characterizations of
non$(\M)$ can be found in [4].
\begin{center}
{\bf References}
\end{center}
\begin{itemize}
\item[{[1]}] T.~Bartoszy\'{n}ski, {\it Combinatorial aspects of measure and
category}, Fund. Math. 127 (1987), pp. 225-239.
\item[{[2]}] T.~Bartoszy\'{n}ski and H.~Judah, {\it Set theory: on the
structure of the real line}, A.~K.~Peters Ltd., Wellesley MA 1995.
\item[{[3]}] A.~W.~Miller, {\it A characterization of the least cardinal
for which the Baire category theorem fails}, Proc. Amer. Math.
Soc. 86 (1982), pp. 498-502. 
\item[{[4]}] A.~Tyszka, {\it On the minimal cardinality of a subset of
$\R$ which is not of first category}, J. Nat. Geom. 17 (2000), pp. 21-28.
\end{itemize}
\begin{flushleft}
{\it Technical Faculty\\
Hugo Ko{\l}{\l}\c{a}taj University\\
Balicka 104, PL-30-149 Krak\'{o}w, Poland\\
rttyszka@cyf-kr.edu.pl\\
http://www.cyf-kr.edu.pl/\symbol{126}rttyszka}
\end{flushleft}
\end{document}